\documentclass[3p]{elsarticle}

\usepackage{amssymb}
\usepackage{natbib}
\usepackage{amsmath}
\usepackage{amsfonts}
\usepackage{graphicx}
\usepackage{mathrsfs}
\usepackage{dcolumn}
\usepackage{bm}
\usepackage{color}
\definecolor{rot}{rgb}{0.75,0.05,0.25}
\definecolor{hellgrau}{gray}{0.5}
\definecolor{blau}{rgb}{0,0,0.7}

\newtheorem{theorem}{Theorem}[section]

\newtheorem{remark}[theorem]{Remark}

\begin{document}

\title{Generalized Thomas-Fermi equations as the Lampariello class of Emden-Fowler equations}

\author[mymainaddress]{Haret C. Rosu\corref{mycorrespondingauthor}}
\cortext[mycorrespondingauthor]{Corresponding author}
\ead{hcr@ipicyt.edu.mx}

\author[mysecondaryaddress]{Stefan C. Mancas}
\ead{mancass@erau.edu}

\address[mymainaddress]{IPICYT,  Instituto  Potosino  de  Investigacion  Cientifica  y  Tecnologica,\\
Camino a la presa San Jos\'e 2055, Col. Lomas 4a Secci\'on,\\ 78216 San Luis Potos\'{\i}, S.L.P., Mexico}
\address[mysecondaryaddress]{Department of Mathematics, Embry-Riddle Aeronautical University,\\ Daytona Beach, FL 32114-3900, USA}

\begin{abstract}
A one-parameter family of Emden-Fowler equations defined by Lampariello's parameter $p$ which, upon using Thomas-Fermi boundary conditions,
turns into a set of generalized Thomas-Fermi equations comprising the standard Thomas-Fermi equation for $p=1$ is studied in this paper.
The entire family is shown to be non integrable by reduction to the corresponding Abel equations whose invariants do
not satisfy a known integrability condition. We also discuss the equivalent dynamical system of equations for
the standard Thomas-Fermi equation and perform its phase-plane analysis. The results of the latter analysis are similar for the whole class.

\end{abstract}
\begin{center}
Physica A 471 (2017) 212-218\\
arXiv: 1607.07072 v4
\end{center}

 \maketitle

\bigskip

\section{Introduction}

The Thomas-Fermi (TF) model \cite{Thomas,Fermi} for the effective electrostatic potential at an arbitrary point in the bulk of a heavy atom emerged as a forerunner of the density functional theory almost ninety years ago, and with some improvements is still in common use for predictions concerning the stability and sizes of heavy atoms, ions, and molecules \cite{LS,Solov}.

The TF equation belongs to the class of Emden-Fowler nonlinear equations which have the form
\begin{equation}\label{ef03b}
\frac{d^2z}{dx^2}=\pm x^{-\lambda-2}z^n~,
\end{equation}
or
\begin{equation}\label{eqef04b}
\frac{d^2y}{dx^2}=\pm x^{-\lambda+n}y^n~,
\end{equation}
where the change of power of the independent variable in the latter equation is obtained via the transformation $z(x)=x ~y\left(\frac 1 x\right)$.
$\lambda$ is a real parameter and $n$ is an integer or fractional parameter. This type of equations became first known in theoretical astrophysics as Lane-Emden and later also as Emden-Fowler equations where they dominated the main stream literature for decades as the basic equations for the Newtonian gravitational potential of a spherically symmetric polytropic gas \cite{Lane,Emden,Fowler}. The astrophysical applications are described generically by the case $\lambda=-2$, with $n$ related to the polytropic index, the dependent variable related to the density of the self-gravitating gas and the independent variable as a dimensionless radius of the gas structure. The minus sign in the right hand side is common for all of the astrophysical applications.
On the other hand, the Thomas-Fermi model \cite{Thomas,Fermi} describing the electronic screening effect in the bulk of a heavy atom leads to the same form of the equation, though for a non-integer $n$ and with the positive sign in the right hand side. In both astrophysics and atomic physics applications, the self-adjoint form of Eq.~(\ref{ef03b}) is used which comes out by means of Kamke's substitutions \cite{Kamke} of the canonical variables  $\eta(\xi)=z(x)$ and $\xi=1/x$
 \begin{equation}\label{EFnormal}
\frac{1}{\xi^2}\frac {d}{d \xi}\left(\xi^2\frac {d\eta}{d\xi}\right)=\pm\xi^{\lambda-2}\eta^{n}~.
 \end{equation}

Other interesting cases correspond to negative $n$'s of the form $n=1-2m$, $\lambda=-m$, in (\ref{ef03b}), with $m$ a positive integer $\geq 2$, when the Emden-Fowler equations can be associated to Ermakov parametric oscillators and their Reid generalization \cite{mr}, while for $m=1$, $\lambda=-2$, and negative sign in the right hand side of (\ref{ef03b}) the equation is known as the `pseudo'-oscillator equation \cite{GL,VG} and has been used to model the path taken by an electron in an electron beam injected into a plasma tube \cite{Act-Sq}.

\medskip

In this paper, we provide a discussion of a class of Emden-Fowler equations presenting fractional powers of both parameters that have been introduced in 1934 by Lampariello \cite{Lamp} which can be considered as {\em generalized Thomas-Fermi equations} since the TF equation is just one particular case of the whole chain and we also maintain the Thomas-Fermi boundary conditions. We show explicitly that these generalized TF equations are related to non integrable Abel equations and therefore no closed solutions are possible for any case of this class. We finally study the autonomous type system of equations to which the Thomas-Fermi equation can be mapped onto, the results of the phase-plane analysis being analogous for all cases.

\section{Generalized Thomas-Fermi equations}

The TF equation is the $p=1$ case of the following class of Emden-Fowler equations introduced by Lampariello \cite{Lamp}
\begin{equation}\label{ef05b}
x^{\frac{p}{p+1}}\frac{d^2y}{dx^2}=y^{\frac{2p+1}{p+1}}~, \qquad {\rm or} \qquad \frac{d^2y}{dx^2}=x^{-\frac{p}{p+1}}y^{1+\frac{p}{p+1}}~,
\end{equation}
with the boundary conditions
\begin{equation}\label{bc}
y(0)=1~, \qquad \lim_{x\rightarrow\infty}y(x)=0~.
\end{equation}
imposed by atomic physics considerations that we keep up for the whole family. An additional condition at the origin is that the physical solution should have finite derivative to the right. This requires $p/(p+1)<1$ \cite{Lamp}, which then provides the interval $(-1,\infty)$ as allowed for the parameter $p$. However, the case $p=0$ is a simple linear equation, while if $p\in (-1,0)$ the origin is not a singular point as it is the case for the standard equations of this class. If one does not impose the additional condition of finite first derivative at the origin, then the interval $p\in (-\infty,-1)$ may be taken into account with some interesting cases such as when the powers of $y$ are positive integers  corresponding to Lane-Emden equations. In the latter cases, we have $1+\frac{p}{p+1}=n\in \mathbb{N}$ which is possible only when $p\in(-2,-1)$.

Comparison of Eq.~\eqref{ef05b} to Eq.~\eqref{eqef04b} implies that $n=2-\frac{1}{p+1}$, and $\lambda =3-\frac{2}{p+1}$.
Thus to obtain the  TF equation, we use $ \lambda=2, n=\frac 3 2$ in Eq.~(\ref{ef03b}) which gives
\begin{equation}\label{ef04}
\frac{d^2z}{dx^2}= x^{-4}z^{\frac 3 2}~
\end{equation}
and according to Eq.~(\ref{EFnormal}) has the self-adjoint form, 
 \begin{equation}\label{EFnormal2}
\frac{1}{\xi^2}\frac {d}{d \xi}\left(\xi^2 \frac{d\eta}{d\xi}\right)=\eta^{\frac 32}~,
 \end{equation}
while with the powers changed takes the standard form
\begin{equation}\label{eq18}
\frac{d^2y}{dx^2}=x^{-\frac 1 2}y^{\frac 3 2}~
\end{equation}
obtained alternatively using $p=1$ in Eq. \eqref{ef05b}.

Equation~\eqref{ef05b} admits a particular solution
\begin{equation}\label{e1}
y_0(x)=\frac{k_p}{x^{1+\frac{2}{p}}}, \qquad k_p=\left[\left(2+\frac 2 p\right)\left(1+\frac 2 p\right)\right]^{1+\frac 1p}~,
\end{equation}
which satisfies the boundary condition \( \displaystyle \lim_{x\rightarrow\infty}y(x)=0\) but not $y(0)=1$.
In addition, Eq.~\eqref{ef05b} is invariant under the rescaling transformations $x \rightarrow ax$ and $y \rightarrow a^{-\left(1+\frac 2 p\right)}y$
for $a\ne 0$ which suggests that there is a function $w(x)=\frac{y(x)}{y_0(x)}$ that satisfies the nonlinear equation
 \begin{equation}\label{eq18bis}
\frac{d^2w}{dx^2}+2\frac{\frac{d y_0}{dx}}{y_0}\frac{dw}{dx}+x^{-\frac {p}{p+1}} y_0^{\frac {p}{p+1} }\left(w-w^{2-\frac {1}{p+1}}\right)=0.
 \end{equation}
 By using $y_0$ from Eq.~\eqref{e1}, we obtain the nonlinear Cauchy-Euler type equation
  \begin{equation}\label{eq190}
x^2\frac{d^2w}{dx^2}-2\left(1+\frac{2}{p}\right)x\frac{dw}{dx}+ \left[\left(2+\frac 2 p\right)\left(1+\frac 2 p\right)\right]\left(w-w^{2-\frac {1}{p+1}}\right)=0.
 \end{equation}
  This equation can be changed into a nonlinear oscillator equation by the rescaled variable $t=\ln x$ to obtain
   \begin{equation}\label{eq191}
\frac{d^2w}{dt^2}-\left(3+\frac{4}{p}\right)\frac{dw}{dt}+k_{p}^{1-\frac{1}{p+1}}w=k_{p}^{1-\frac{1}{p+1}}w^{2-\frac {1}{p+1}}~,
 \end{equation}
where we have moved the nonlinearity on the right hand side. These nonlinear oscillators have a negative constant `damping' ratio
$\zeta_p=-\left(\frac{3}{2}+\frac{2}{p}\right)k_{p}^{\frac{1}{2}\left(\frac{1}{p+1}-1\right)}$ and a specific stiffness of $\kappa_p=k_{p}^{1-\frac{1}{p+1}}$ of the same value as the specific strength of the nonlinearity. The interesting feature of these nonlinear oscillators is that when $p$ is increased from $1$ (the TF case) to $\infty$, the `damping' ratio varies only from $\zeta_1=-7\sqrt{ 3}/12\approx-1.01$ to $\zeta_\infty=-3 \sqrt{2}/4\approx -1.06$, whereas the specific stiffness decreases six times, from $\kappa_1=12$ to $\kappa_\infty=2$, in parallel  with the same decrease of the specific strength of the nonlinearity while the latter raises its power index from three halves to quadratic.
To display the limits of variation of the coefficients of this class of nonlinear equations, we use $p=1$ and $p=\infty$ in  Eqs.~\eqref{eq190} and
\eqref{eq191} that gives
 \begin{equation}\label{eq18ab}
x^2\frac{d^2w}{dx^2}-6 x \frac{dw}{dx}+12 w=12 w^{\frac 32}~, \qquad  x^2\frac{d^2w}{dx^2}-2 x \frac{dw}{dx}+2 w=2w^{2}~,
 \end{equation}
which in the rescaled variable corresponds to
 \begin{equation}\label{eq19}
\frac{d^2w}{dt^2} -7 \frac{dw}{dt}+12 w=12w^{\frac 32}~, \qquad \frac{d^2w}{dt^2} -3\frac{dw}{dt}+2 w=2w^{2}~,
 \end{equation}
 respectively.
We also notice that the characteristic equation of the corresponding linear oscillator has the roots $r_1=2+\frac 2 p$ and
$r_2=1+\frac 2 p$. In the case of the TF equation, $r_1=4$ and $r_2=3$; moreover, their squared product is $(r_1r_2)^2= k_1=144$.

\subsection{Perturbation of the particular solution}

 Now, let us perturb the particular solution by a function $\epsilon(x)=y(x)-y_0(x)$, assumed infinitesimal. By substituting into Eq. \eqref{ef05b} and expanding with respect to $\epsilon$ we obtain

  \begin{equation}\label{eps0}
x^{\frac{p}{p+1}}\frac{d^2\epsilon}{dx^2}=\frac{2p+1}{p+1}y_0^{\frac{p}{p+1}}\epsilon+\frac{p(2p+1)}{2(p+1)^2}y_0^{-\frac{1}{p+1}}\epsilon^2-\frac{p(2p+1)}{6(p+1)^3}y_0^{-\left(1+\frac{1}{p+1}\right)}\epsilon^3+{\cal O}(\epsilon^4)~.
 \end{equation}
 Using the particular solution $y_0$ the above reduces to
   \begin{equation}\label{eps1}
\frac{d^2\epsilon}{dx^2}=\frac{2p+1}{p+1}k_p^{\frac{p}{p+1}}\frac{1}{x^2}\epsilon+\frac{p(2p+1)}{2(p+1)^2}k_p^{-\frac{1}{p+1}}x^{\frac 2 p -1}\epsilon^2-\frac{p(2p+1)}{6(p+1)^3}k_p^{-\left(1+\frac{1}{p+1}\right)}x^{\frac 4p}\epsilon^3+{\cal O}(\epsilon^4).
 \end{equation}
If we keep now only the linear term,  we obtain the Cauchy-Euler equation
\begin{equation}\label{eps10}
\frac{d^2\epsilon}{dx^2}=\frac{2p+1}{p+1}k_p^{\frac{p}{p+1}}\frac{1}{x^2}\epsilon
 \end{equation} from which we get
   \begin{equation}\label{eps10b}
\frac{d^2\epsilon}{dx^2}=\frac{18}{x^2}\epsilon
 \end{equation}
 as a particular case when $p=1$ with general solution
  \begin{equation}\label{eps3}
 \epsilon(x)=c_1 x^{\frac 1 2(1-\sqrt{73})}+c_2x^{\frac 1 2(1+\sqrt{73})}\approx c_1x^{-3.772}+c_2x^{4.772}~.
 \end{equation}
Lampariello used the negative power term in (\ref{eps3}) to argue on the asymptotic behavior of solutions at infinity. Since it is the only one which is infinitesimal at infinity and moreover since its power exponent is more negative than the power exponent of $y_0$, which is valid for any $p$,
then any solution $y(x)$ which is infinitesimal at infinity differs from $y_0$ by an infinitesimal of higher order with respect to $y_0$. One concludes that the solution $y(x)$ itself should differ from $y_0$ by an infinitesimal of order higher to the order of $y_0$
and then one may assume that the quotient $w(x)=\frac{y(x)}{y_0(x)} \rightarrow 1$ for $x\rightarrow\infty$.
On the other hand, the presence of the irregular positive power term is indicative of possible perturbative instabilities and may explain the nonexistence of more extended Thomas-Fermi structures. In addition, we will see in the last section that the same powers as in (\ref{eps3}) correspond to the eigenvalues of the Jacobian matrix at the saddle point.

\section{Abel equation as an ingredient for generalized Thomas-Fermi equations}

To write Eq.~(\ref{eq191}) as an Abel equation of the second kind
 we proceed as in \cite{Ve} and let $\frac {dw}{dt}=s(w)$. This leads to
    \begin{equation}\label{eq1911}
    s\frac{ds}{dw}-\left(3+\frac{4}{p}\right)s+ \left[2\left(1+\frac 1 p\right)\left(1+\frac 2 p\right)\right]\left(w-w^{2-\frac {1}{p+1}}\right)=0
 \end{equation}
which corresponds to
 \begin{equation}\label{eq19a}
s\frac{ds}{dw}-7s +12 w(1-\sqrt w)=0~
 \end{equation}
for the TF equation when $p=1$.

\subsection{Abel Invariant}

We use the inverse transformation $s(w)=\frac{1}{z(w)}$ to obtain the Abel equation of the first kind
 \begin{equation}\label{eq20}
\frac{dz}{dw}=f_2 z^2+f_3 z^3=-  \left(3+\frac{4}{p}\right) z^2+\left[2\left(1+\frac 1 p\right)\left(1+\frac 2 p\right)\right]\left(w-w^{2-\frac {1}{p+1}}\right)z^3~.
 \end{equation}

 This equation has the invariant
  \begin{equation}\label{eq21}
\Phi_p(w)=\frac 1 3 \left(f_3\frac{d f_2}{d w}-f_2\frac{d f_3}{d w}\right)+\frac{2}{27}f_2^3= \frac{2(3p+4)}{27 p^3}\left[(3p+2)-9(p+2)(2p+1)w^{1-\frac{1}{p+1}}\right]~.
 \end{equation}

\begin{remark}
Abel equations are integrable if the invariant satisfies the integrability condition
   \begin{equation}\label{eq22}
f_3\frac{d \Phi_p}{d w}+\left(f_2^2-3\frac{d f_3}{d w}\right)\Phi_p=3 \alpha \Phi_p^{\frac 53},
 \end{equation}
for some constant $\alpha$.
In that case, the solution is
\begin{equation}\label{sol1}
z=\frac{3 u {\Phi_p}^{\frac 1 3}-f_2}{3 f_3}
\end{equation}
where  $u$ is found by one quadrature
 \begin{equation}\label{sol2}
\int \frac {du}{u^3-\alpha u+1}+C=\int \frac{{\Phi_p}^{\frac 2 3}}{f_3}dw~.
\end{equation}
\end{remark}
However, the structure of the invariant $\Phi_p$ in Eq.~(\ref{eq21}) leads to an integrability condition of the form
 \begin{equation}\label{intcond}
 c_1 w^q + c_2 w^{2q} +
  c_3 \left(c_4 - c_5 w^q\right)^{\frac 53} \alpha = c_6~,
 \end{equation}
where the constants $c_i=c_i(p)$ and $q=\frac{p}{p+1}$, and there is no constant $\alpha$ that satisfies (\ref{intcond}).
In particular, the Abel equation corresponding to the TF equation has the invariant
   \begin{equation}\label{eq201}
\Phi_1(w)=\frac{70}{27} - 42 \sqrt w
 \end{equation}
 and the corresponding integrability condition
 \begin{equation}\label{eq23}
 3807 \sqrt{w} + 11664 w +
  14^{\frac 2 3} \left(5 - 81 \sqrt{w}\right)^{\frac 53} \alpha = 195
   \end{equation}
   is not satisfied by any constant $\alpha$.

Therefore, we will instead use a special transformation due to Lampariello \cite{Lamp} which will change the Thomas-Fermi equation into a first-order equation obtained by Majorana \cite{Maj}.

\subsection{The Lampariello transformation}

 Let us rewrite Eq. \eqref{eq1911} as
     \begin{equation}\label{eq1912}
    s\frac{ds}{dw}=\left(3+\frac{4}{p}\right)s+ \left[2\left(1+\frac 1 p\right)\left(1+\frac 2 p\right)\right]w \left(w^{\frac {p}{p+1}}-w\right)
 \end{equation}
 and use  $w(\tau)=\tau^{(p+1)(p+2)}$ with
     \begin{equation}\label{eq24}
s(w)=\left(1+\frac{2}{p}\right) \left[1-\tau^{p(p+1)}u(\tau)\right]w
   \end{equation}
   in Eq.~(\ref{eq1912}). Then we obtain
       \begin{equation}\label{eq250}
\frac{d u}{d\tau}=-2(p+1)^2~\frac{\tau^{p-1}\left(1-\tau^{p^2} u^2\right)}{1-\tau^{p(p+1)}u}~,
   \end{equation}
   which in the TF case takes the form
   \begin{equation}\label{eq25}
   \frac{d u}{d\tau}=\frac{-8(1-\tau u^2)}{1-\tau^2u}~.
   \end{equation}
  This equation is the same as Eq.~(27) in Esposito's paper on the results of Majorana \cite{Maj} and it is an important intermediate step in Majorana's approach. A power expansion of the function $u$ in the variable $(1-\tau)$ led Majorana to his parametric solution of the TF equation adapted to the phenomenological requirements of atomic physics.

\section{An autonomous two-dimensional ODE system of the TF equation}

This section deals essentially with the TF equation since the results are similar for the other nonunity values of $p$.
By transforming the self-adjoint form of the TF equation into a dynamical system, one can classify the solutions based on linear stability analysis.
This mapping can be achieved by using the transformations on Eq. \eqref{EFnormal} given by Jordan and Smith in their book \cite{js} 
\begin{eqnarray}\label{l1}
\begin{array}{ll}
&X=\xi ~\frac{\eta_{\xi}}{\eta}\label{v1}\\
&Y=\xi^{\lambda-1} ~\frac{\eta^n}{\eta_{\xi}}~,\label{v2}
\end{array}
\end{eqnarray}
with  $\xi=e^t$, which will turn the self-adjoint form of TF Eq.~(\ref{EFnormal2}) into the autonomous two-dimensional ODE system
\begin{eqnarray}\label{l2}
\begin{array}{ll}
& \frac{dX}{dt}=-X(1+X-Y)=M(X,Y)\\
&\frac{dY}{dt}=Y(\lambda+1+n X- Y)=N(X,Y)~,
\end{array}
\end{eqnarray}
with four equilibrium points given by
$$
\left\{ (X_0,Y_0)=(0,0);  ~(X_1,Y_1)=(-1,0);~  (X_2,Y_2)=(0,\lambda+1); ~(X_3,Y_3)=\left(-\frac{\lambda}{n-1},1-\frac{\lambda}{n-1}\right) \right\}.
$$
We mention that a discussion of a similar system of equations obtained with the same transformations applied to the non self-adjoint form of the TF equation has been provided by Hille a long time ago \cite{H70}. However, his discussion is focused on adding more intuition to the series expansions of the TF solutions for large and small values of the independent variable rather than on the phase-plane analysis as we do next.

Following standard methods of phase-plane analysis, we  use the linear approximation of the equilibrium points to classify them. The Jacobian matrix of (\ref{l2}) is

\begin{eqnarray}\label{17}
J=\left[\begin{array}{cc}
\frac{\partial M}{\partial X}& \frac{\partial M}{\partial Y}\\
\frac{\partial N}{\partial X}& \frac{\partial N}{\partial Y}\\
\end{array}\right]=\left[\begin{array}{cc}
-1-2X+ Y& X\\
nY& \lambda+1+n X- 2Y\\
\end{array}\right]
\end{eqnarray}
and the characteristic polynomial of the Jacobian matrix is
\begin{equation} \label{charac1}
j(\theta)=\theta^2-\delta_1 \theta+\delta_2=0~.
\end{equation}

\begin{figure}
\centering
\includegraphics[width=0.75\textwidth]{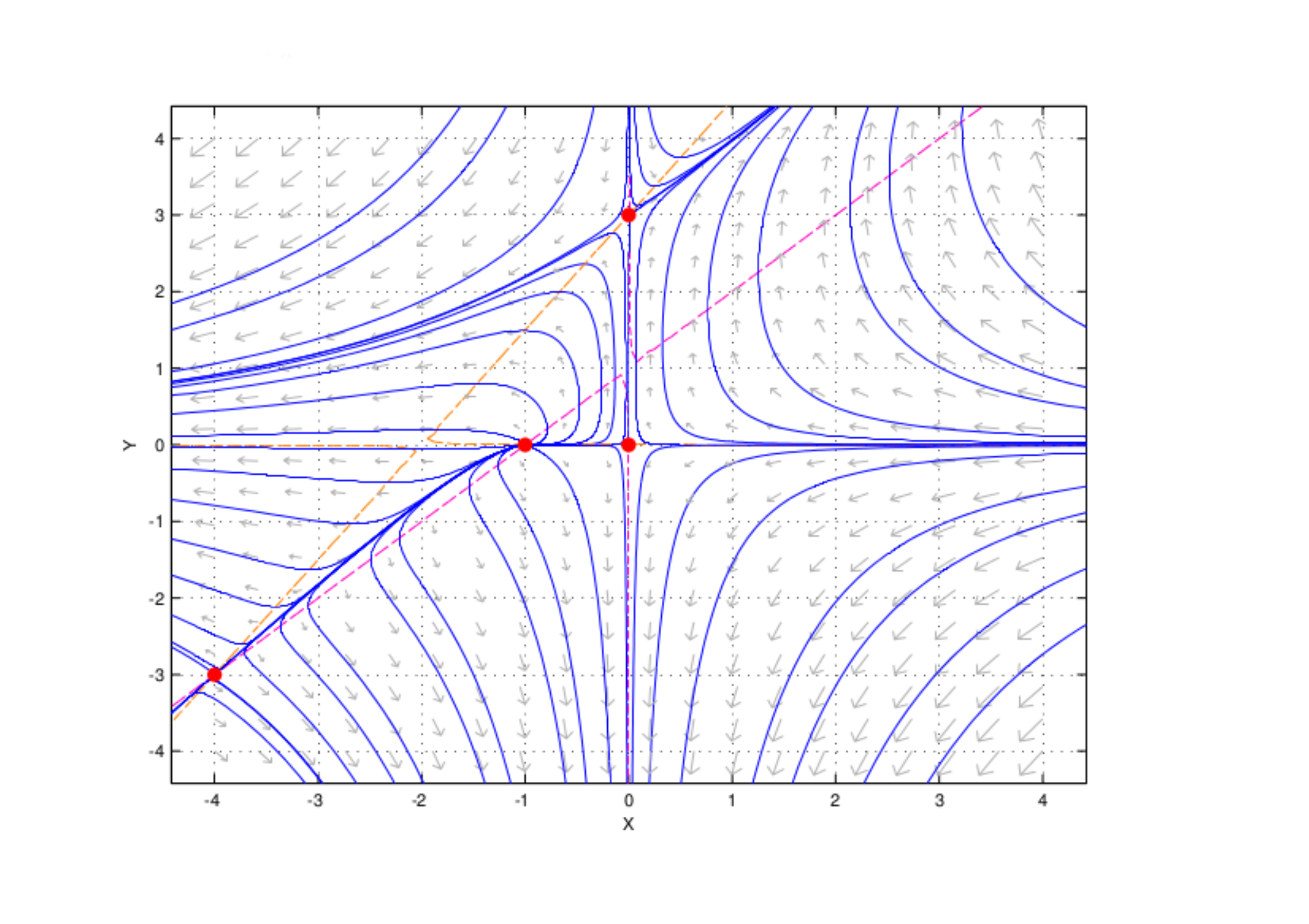}
\caption{\small{Phase-plane portrait of the autonomous ODE system (\ref{l2}) for the TF equation for $n=\frac 3 2$ and $\lambda=2$.}}
\label{figure1}
\end{figure}

The equilibrium points will be classified according to signs of the trace $\mathrm{Tr}(J)=\delta_1=\frac{\partial M}{\partial X}+ \frac{\partial N}{\partial Y}$, determinant $\mathrm{Det}(J)=\delta_2=\frac{\partial M}{\partial X}\frac{\partial N}{\partial Y}-\frac{\partial M}{\partial Y}\frac{\partial N}{\partial X}$, and discriminant $\Delta=\delta_1^2-4\delta_2$, all evaluated at $(X_i,Y_i)$.
For TF equation, $n=\frac 32, \lambda=2$, the results are presented in Table \ref{Tab1} and the phase-plane portraits in Fig.~\ref{figure1}.

\begin{table}[htb!] %
\begin{center}
\begin{tabular}{|c|c|c|c|c|c|c|}
\hline Fixed Points & $\delta_1$ & $\delta_2$ & $\Delta$ & Type \\
\hline
\hline
$(X_0,Y_0)$ & $2$ & $-3$ & $16$   & saddle\\
\hline
$(X_1,Y_1)$  & $\frac 52 $ & $\frac 3 2$ & $\frac 1 4$  &  stable node\\
\hline
$(X_2,Y_2)$ & -1 & $-6$  & $25$ & saddle \\
\hline
$(X_3,Y_3)$ & $7$ & $-6$ & $73$ & saddle\\
\hline
\end{tabular}
\end{center}
\caption{General equilibrium points of the predator-prey system (\ref{l2}).}
\label{Tab1}
\end{table}

Now let us consider the only nontrivial fixed point $(X_3,Y_3)=(-4,-3)$. Since in the general case $(X_3(p),Y_3(p))=(-3-1/p,-2-1/p)$, this fixed point moves towards $(-3,-2)$ when $p\rightarrow \infty$. Using the transformation given by \eqref{l1} we have $XY=\xi^\lambda \eta^{n-1}\Rightarrow X_3Y_3=12=\xi^2 \eta^{\frac 1 2}=\frac{\sqrt{z(x)}}{x^2}$. Since $z(x)=xy(\frac 1 x)$, the particular solution $y_0(x)=144/x^3$ can also be obtained using the saddle point.
Also, the Jacobian at the fixed point is

\begin{eqnarray}\label{170}
J_{(X_3,Y_3)}=\left[\begin{array}{cc}
4& -4\\
-\frac 9 2 & 3\\
\end{array}\right]
\end{eqnarray}
with characteristic polynomial $j_{(X_3,Y_3)}=\theta^2-7 \theta-6=0$. The eigenvalues
are  $\theta_{1,2}=\frac{7\pm\sqrt{73}}{2}=3+\theta_{0_{1,2}}$, where the  perturbed values from the integer part $3$ are $\theta_{0_{1,2}}=\frac{1\pm\sqrt{73}}{2}$. The integer part of the eigenvalue, $3$,
corresponds to the power of the particular solution $y_0$, while the perturbed eigenvalues correspond to the powers of the perturbed solutions $\epsilon$ from Eq.~\eqref{eps3}. For the first eigenvalue $\theta_1$ the flow is in the direction of
$
\vec u_1=\left[\begin{array}{c}
1\\
1-\frac {\theta_1}{4} \\
\end{array}\right]\approx \left[\begin{array}{c}
1\\
-0.943 \\
\end{array}\right]
$,
while for  the second eigenvalue $\theta_2$ the flow is in the direction of
$
\vec u_2=\left[\begin{array}{c}
2(3-\frac{\theta_2}{9})\\
1 \\
\end{array}\right]\approx\left[\begin{array}{c}
6.171\\
1 \\
\end{array}\right]~.
$
\bigskip

\section{Conclusion}

The generalized class of TF equations labeled by the Lampariello parameter $p\geq 0$, with the TF equation corresponding to the case $p=1$, has been shown to be nonintegrable by reduction to the corresponding class of Abel equations.
A different reduction of the TF equation to an autonomous system of first-order equations shows that the system possesses a single stable node and three unstable saddle points in the phase plane, a result which is valid for the whole class of generalized TF equations.

We surmise that the TF equations in this class may be used to describe quantum systems which display intrinsic deviations from Fermi or Bose statistics \cite{Greenberg} or in some astrophysical or cosmological context in which the quantum electrostatics may present more interweaved screening effects.

Finally, we have also noticed that if the condition of boundedness of the first derivative of the solution at the origin is relaxed then one can consider this family of equations in the direction $p<-1$ with direct connection to the Lane-Emden equations in astrophysics only for $p\in(-2,-1)$.

\bigskip

\noindent {\bf Acknowledgment}: The authors thank the referees for comments which helped to perform a substantial improvement of this work.


\end{document}